\documentclass[12pt]{article}

\usepackage{diagbox}
\usepackage{mathtools}
\usepackage{bbm}
\usepackage{latexsym}
\usepackage{epsfig}
\usepackage{amsmath,amsthm,amssymb,enumerate}

\usepackage[a-1b]{pdfx}
\usepackage{hyperref}

\usepackage{color}

\def\cH{{\mathcal H}}
\def\nn{\nonumber}
\def\a{\alpha}   
\def\e{\varepsilon}    \def\g{\gamma}

 \def\m{\mu}

\newtheorem{theorem}{Theorem}

\newtheorem{claim}{Claim}

\newcommand{\brac}[1]{\left(#1\right)}

\newcommand{\bfrac}[2]{\left(\frac{#1}{#2}\right)}

\def\cE{{\cal E}}

\def\E{\mathbb{E}}

\def\Pr{\mathbb{P}}

\newcommand{\ignore}[1]{}

\def\cE{{\mathcal E}}

\def\cH{{\mathcal H}}

\def\nn{\nonumber}

\def\tm{\widetilde{m}}

\usepackage{tikz}
\usetikzlibrary{decorations.pathmorphing}
\usetikzlibrary{positioning}
\usetikzlibrary{arrows,automata}
\usetikzlibrary{shapes.misc}
\usetikzlibrary{backgrounds}
\usetikzlibrary{arrows,shapes}

\begin{document}
\author{Patrick Bennett\thanks{Research supported in part by Simons Foundation Grant \#426894.}\and Alan Frieze\thanks{Research supported in part by NSF grant DMS1952285}}

\date{}
\title{Some online Maker-Breaker games}
\maketitle

\begin{abstract}
    We consider some Maker-Breaker games of the following flavor. We have some set $V$ of items for purchase. Maker's goal is to purchase some member of a given family $\cH$ of subsets of $V$ as cheaply as possible and Breaker's goal is to make the purchase as expensive as possible. Each player has a pointer and during a player's turn their pointer moves through the items in the order of the permutation until the player decides to take one. We mostly focus on the case where the permutation is random and unknown to the players (it is revealed by the players as their pointers move). 
\end{abstract}

\section{Introduction}

In this note, we consider a few related games with two players, Maker and Breaker. In general, we have some set $V$ of items for purchase, and Maker's goal is to purchase some member of a given family $\cH$ of subsets of $V$ (as a subset of the items that Maker purchases). We will focus on cases where Maker can guarantee (or ``almost'' guarantee) that they are able to get some member of $\cH$. But each item has a cost, and Maker would also like to minimize the total cost of the items they purchase. Breaker's goal will be to maximize the cost Maker needs to pay. In this sense Breaker does not play the same role as the type of player usually called ``Breaker,'' whose goal is to entirely prevent Maker from collecting some set of items. Instead, in our games Breaker is resigned to allowing Maker to succeed in collecting a desired set of items, but wants to make it expensive.  

The items in $V$ will be ordered in a uniform random permutation. Each item has a uniform random cost in the interval $[0, 1]$. The permutation of $V$ and the costs are unknown to the players in the beginning. Each player has a pointer which in the beginning points at the first item of $V$ in the permutation. As the game progresses the pointers will move forward through the permutation. Once either player's pointer moves past an item of $V$, both players know the item's cost and what position it is in the random permutation of $V$. Of course, our analysis of these games will often involve ``revealing'' the random costs and the permutation as the pointers move forward, leaving the rest of the permutation and random costs ``unrevealed.'' The players alternate turns and Breaker's turn is first. Breaker's pointer moves forward one item at a time, and as Breaker's pointer moves over an item they see its cost and decides to either take it or reject it. Once Breaker has taken a total of $b$ items (or reached the end of $V$), it is Maker's turn. Maker's pointer is still at the first item, but Maker knows everything Breaker has revealed about the permutation and costs up to wherever Breaker's pointer is now. Maker examines the items in order and for each one that hasn't already been taken by Breaker, Maker can purchase it or reject it. Thus Maker can either take an item rejected by Breaker, or Maker could wait until Maker's pointer moves past Breaker's to start revealing new items until Maker decides to purchase one (or reaches the end of $V$). Then Breaker gets another turn, starting wherever Breaker's pointer was when Breaker's last turn ended. And so on.

In the first game we consider, $\cH$ is just the set of all singletons, i.e. Maker's goal is to purchase one item. We call this game \textsc{item}~since Maker just wants one item. We prove the following. 

\begin{theorem}\label{thm:item}
    For the game \textsc{item}~with $b=1$  there exists a strategy for Maker which always purchases an item and where the expected cost paid by Maker is at most $(4+o(1))/n$. For all     $1 \le b \le n / \log^4 n$ there exists such a strategy for Maker  where the expected cost is at most $O\brac{\tfrac{b \log ^2 b }{ n}}.$ Furthermore, there exists a strategy for Breaker for which Maker's expected cost is $ \Omega\brac{\tfrac{b  }{ n}}.$
\end{theorem}

Sometimes we will make an additional restriction on Breaker to prevent them from getting too far ahead of Maker. Maker's strategy will often involve using the early items in the permutation to build a collection that, for example, contains all but one item of $e$ for many sets $e \in \cH$. Maker will then look in the later items in the permutation to find a single low-cost item which would complete such a set $e$. But if Breaker knows that Maker plans to follow such a strategy then Breaker can pass over the early part of the permutation, observing the structure that Maker will build from it, and spoil Maker's strategy in the later parts of the permutation by removing the items Maker will want. In some situations it seems Maker is powerless to stop Breaker from doing this, since Maker will have to spend a significant amount of turns in the earlier part of the permutation in order to purchase the necessary items. While Maker is stuck in the earlier part of the permutation, Breaker is free to cause this kind of trouble in the later part. Thus, we sometimes add a restriction to prevent Breaker from doing this. We say the game is $p$-phase restricted if the permutation of items is divided into $k$ contiguous pieces we call the {\em phases}, each consisting of $\lfloor n/p \rfloor$ or $\lceil n/p \rceil$ steps, such that Breaker is not allowed to start a phase until Maker has reached the end of the previous phase. 

We now consider the case of trying to purchase a subgraph. More precisely, our set of items $V$ will be the edge set of $K_n$, and $\cH$ some collection of sets of edges in $K_n$. Our next result is when $\cH$ is the collection of triples of edges that form a triangle in $K_n$ (more generally, $\binom k2$-tuples of edges forming a $k$-clique). The game where Maker tries to build a $k$-clique will be called \textsc{$k$-clique}.

\begin{theorem}\label{thm:clique}
For the game \textsc{$3$-clique} with $b \le n^{2/3}/10$ and no phase restriction, there is a strategy for Maker which w.h.p.~succeeds in purchasing a $3$-clique and spends at most $O(b n^{-1/3} \log n)$. For each fixed $k\ge 3$ and for all $b = o(n^{11\a_k/4}),\;\a_k=\tfrac{1}{11 \cdot 2^{k-5}-1}$ there exists a $k$-phase restricted strategy for Maker which w.h.p.~succeeds in purchasing a $k$-clique and spends at most $O\brac{b n^{-\a_k} \log ^2 b}$.
\end{theorem}

Next we consider the problem where Maker wants a path from $u$ to $v$ for two given vertices $u, v$. For the case $b=0$, Frieze and Pegden \cite{FP18} proved that this can be done with an expected cost of $n^{-2/3 + o(1)}$. 

\begin{theorem}\label{thm:paths}
   For $b \le n^{1-1/\log \log n}$ there is an $O(\log \log n)$-phase restricted strategy for Maker to purchase a path between two given vertices which succeeds with high probability with a total cost of at most $bn^{-2/3+o(1)}$
\end{theorem}

Finally, we consider an ``order-restricted'' version of the box game studied by Chv\'{a}tal and Erd\H{o}s \cite{CE78}. The original box game from \cite{CE78} is also a Maker-Breaker game (although in our description we will swap the names of Maker and Breaker compared to the description in \cite{CE78}). We start with $n$ boxes, each having $m$ balls. The players take turns, Maker going first. Maker takes one ball each turn, and Breaker takes $b$ balls. Maker's goal is to take one ball from every box. Chv\'{a}tal and Erd\H{o}s \cite{CE78} gave the exact value $m_0=m_0(n, b)$ such that Maker has a winning strategy if and only if $m \ge m_0$. For large $n$, this $m_0$ is about $b \log n$. 

In the ``order-restricted'' version we consider here, the balls have some ordering which is initially unknown to the players (the ordering is revealed as the game is played). On each player's turn they are offered the balls according to the order, picking up from wherever they left off at their last turn (and of course they cannot take a ball that was already taken by the other player). We call our version \textsc{BoxGame}. We consider the case where the balls are randomly ordered as well as the case of general orderings.

\begin{theorem}\label{thm:boxgame}
If the ordering of balls is arbitrary (possibly adversarial favoring Breaker, and adaptive to the moves made so far), then for all $b \ge 1$, Maker has a winning strategy for \textsc{BoxGame} if and only if $m \ge bn+1$. If the ordering is random and $b \ge b_0=100\e^{-2}\log n$, then for any fixed $\e >0$ and $m \le (1-\e)bn$, there is a strategy for Breaker to win with high probability. 
\end{theorem}
    
\section{Purchasing one item}

In this section we prove Theorem \ref{thm:item}. The lower bound is easy so we do it first. 

\subsection{Lower bound}
Breaker's strategy is to take any item that costs at most $b/2n$. The expected number of such items is $b/2$, so by Markov's inequality there are at most $b$ of them with probability at least $1/2$. In that event Breaker will get all of those items and then Maker has to pay at least $b/2$. Thus Maker's expected cost is at least $b/4$. This completes the proof of the lower bound. 

 \subsection{Upper bound for $2 \le b \le n / \log^4 n$} We handle the case where $2 \le b \le n / \log^4 n$. Let
\begin{equation}\label{eqn:alphadef}
    \a:= 10+10 \lceil \log b \rceil
\end{equation}
Set 
\begin{equation}\label{eqn:nmidef}
    N:= \frac{n}{b+1}, \qquad  m_{i} := \frac{\a}{N+\a - i}.
\end{equation}
Maker's strategy is as follows. First we break the $n$ steps into $b+1$ {\em phases} of $N$ consecutive steps. In each phase, Maker will try to take the first item $i$ such that $c_i \le m_i$ (where the value of $i$ restarts at each phase so each phase has $i=1, \ldots, N$). Since $m_N=1$ it is guaranteed that each phase has such a value $i$. Of course when we say Maker ``tries'' to take the item, we mean that Maker takes it unless Breaker already has. But, since there are $b+1$ phases, Breaker must allow Maker to have one of them.   

For $1 \le i \le N$, let $\cE_{i}$ be the event that out of all $b+1$ phases, $i$ is the largest indexed item that Maker has tried to take. In other words, in one of the phases Maker has tried to take item $i$ and in all the other phases Maker has tried to take some item at most $i$. In the event $\cE_{i}$ Maker gets some item indexed at most $i$ which could cost at most $m_i$. The expected cost paid by Maker is at most
\begin{equation}\label{eqn:mean1item}
        \sum_{i=1}^N \Pr\brac{\cE_{i}} m_{i}.
\end{equation}
We now start trying to bound \eqref{eqn:mean1item} from above. We have
\begin{equation}\label{eqn:ei}
     \Pr\brac{\cE_{i}} \le  (b+1) \left[ \brac{\prod_{j=1}^{i-1} (1-m_{j})} m_{i}\right] \left[ 1- \brac{\prod_{j=1}^{i} (1-m_{j})}\right]^b.
\end{equation}
We address the large values of $i$ as follows. Set $x:= \a^2 \log N$, and note that since we have $1 \le b \le n / \log^4 n$, \eqref{eqn:alphadef} and \eqref{eqn:nmidef} it follows that 
\begin{equation}\label{eqn:xasymptotics}
    \Omega(1)  \le \a \le O(\log n), \qquad N = \Omega(\log^4 n), \qquad \a = o(N), \qquad \a^2 = o(x).
\end{equation}
For $i > N-x$ we have
\begin{equation}\label{eqn:PEbound}
    \Pr\brac{\cE_{i}} \le  (b+1)  \prod_{j=1}^{N-x} (1-m_{j})  \le (b+1)\exp \left\{ -\sum_{j=1}^{N-x} m_{j} \right\}
\end{equation}
Now by the approximation of the harmonic series 
\begin{equation}\label{eqn:harm}
    \sum_{\ell = 1}^m \frac 1\ell = \log(m) + \g + O\brac{\frac{1}{m}}
\end{equation}
 we have
\begin{align}
  \sum_{j=1}^{N-x} m_{j} &= \sum_{j=1}^{N-x} \frac{\a}{N+\a - j}\nonumber\\
  & = \a \bigg( \log(N+\a-1) - \log(x+\a-1) + O\brac{\frac{1}{x+\a-1}}\bigg)\nonumber\\
  &= - \log \brac{\frac{x+\a-1}{N+\a-1}}^\a+ o(1)\label{eqn:summbigi}
\end{align}
where on the last line we have used $\frac{\a}{x+\a-1} = o(1)$ which follows from \eqref{eqn:xasymptotics}. Thus the total contribution to \eqref{eqn:mean1item} by all terms $i > N-x$ is at most $x$ times \eqref{eqn:PEbound}, which in light of \eqref{eqn:summbigi} is at most
\begin{equation}\label{eqn:totalbigi}
   O\bigg( x (b+1) \brac{\frac{x+\a-1}{N+\a-1}}^\a \bigg)=  O\bigg( x (b+1) \brac{\frac{2x}{N}}^\a \bigg).
\end{equation}
Now we address the rest of the values $i$. For $i \le N-x$ we have for all $j \le i$ that $m_j \le m_{N-x} = \a / (x+\a) = o(1)$. Thus, using the Taylor series for $\log(1-x)$ gives
\begin{align}
  \prod_{j=1}^{i} (1-m_{j}) & = \exp \left\{ \sum_{j=1}^{i} \log(1-m_{j})\right\}   = \exp \left\{ -\sum_{j=1}^{i} (m_{j} + O(m_j^2))\right\}.\label{eqn:est1}
\end{align}
Using \eqref{eqn:harm} again we have
\begin{align}
  \sum_{j=1}^{i} m_{j} &= \sum_{j=1}^{i} \frac{\a}{N+\a - j}\nonumber\\
  & = \a \bigg( \log(N+\a-1) - \log(N+\a-i-1) + O\brac{\frac{1}{N+\a-i-1}}\bigg)\nonumber\\
  & =  -\log\brac{\frac{N+\a-i-1}{N+\a-1}}^\a + o(1)\label{eqn:summsmalli}
\end{align}
where on the last line we used $\frac{\a}{N+\a-i-1} \le \frac{\a}{x+\a-1} = o(1)$ by \eqref{eqn:xasymptotics}. 
Using
\[
\sum_{\ell = 1}^m \frac 1{\ell^2} = \frac{\pi^2}{6} + O\brac{\frac{1}{m}}
\]
we have
\begin{align}
     \sum_{j=1}^{i} m_{j}^2 &= \sum_{j=1}^{i} \frac{\a^2}{(N+\a - j)^2} = O\brac{\frac{\a^2}{N+\a-i}}=o(1)\label{eqn:summsquared}
\end{align}
since $\frac{\a^2}{N+\a-i-1} \le \frac{\a^2}{x+\a-1} = o(1)$ by \eqref{eqn:xasymptotics}.
Now by \eqref{eqn:summsmalli} and \eqref{eqn:summsquared}, \eqref{eqn:est1} is 
\[
\exp \left\{\log\brac{\frac{N+\a-i-1}{N+\a-1}}^\a + o(1)\right\}= (1+o(1)) \brac{\frac{N+\a-i-1}{N+\a-1}}^\a .
\]
Therefore for $i \le N-x$, \eqref{eqn:ei} is 
\begin{align*}
    &(1+o(1)) (b+1) m_i \brac{\frac{N+\a-i-1}{N+\a-1}}^\a  \left[ 1- (1+o(1))\brac{\frac{N+\a-i-1}{N+\a-1}}^\a \right]^b\\
    & \le 2(b+1) m_i \brac{\frac{N+\a-i-1}{N+\a-1}}^\a  \left[ 1- \frac12 \brac{\frac{N+\a-i-1}{N+\a-1}}^\a \right]^b\\
    & \le 2(b+1) m_i \brac{\frac{N+\a-i-1}{N+\a-1}}^\a  \exp \left\{ - \frac12 b\brac{\frac{N+\a-i-1}{N+\a-1}}^\a \right\}
\end{align*}
Thus, the contribution to \eqref{eqn:mean1item} by terms $i \le N-x$ is at most 
\begin{align}
    &2(b+1) \sum_{i=1}^{N-x}  m_i^2 \brac{\frac{N+\a-i-1}{N+\a-1}}^\a   \exp \left\{ - \frac12 b\brac{\frac{N+\a-i-1}{N+\a-1}}^\a \right\}\nonumber\\
    & = 2(b+1) \sum_{i=1}^{N-x}  \brac{\frac{\a}{N+\a-i}}^2 \brac{\frac{N+\a-i-1}{N+\a-1}}^\a   \exp \left\{ - \frac12 b\brac{\frac{N+\a-i-1}{N+\a-1}}^\a \right\}\nonumber\\
    & \le \frac{2\a^2(b+1)}{(N+\a-1)^2}\sum_{i=1}^{N-x}   \brac{\frac{N+\a-i-1}{N+\a-1}}^{\a-2}   \exp \left\{ - \frac12 b\brac{\frac{N+\a-i-1}{N+\a-1}}^\a \right\}\label{eqn:smalli}
\end{align}
Letting 
\[
f(t):= t^{\a-2} e^{-\frac12 bt^\a},
\]
elementary calculus implies the maximum of $f$ over $0 \le t \le 1$ occurs when $t=\brac{\frac{2(\a-2)}{\a b}}^{1/\a}$. Thus each term of \eqref{eqn:smalli} is at most \[
f\brac{\brac{\frac{2(\a-2)}{\a b}}^{1/\a}}= \brac{\frac{2(\a-2)}{\a b}}^{1-2/\a} e^{-1+2/\a}  =O\brac{ b^{-1+2/\a} }= O(b^{-1})
\]
where we have used \eqref{eqn:alphadef}. 
Thus \eqref{eqn:smalli} is at most
\begin{equation}\label{eqn:smalli2}
      \frac{2\a^2(b+1)}{(N+\a-1)^2} \cdot (N-x) \cdot  O(b^{-1}) = O\brac{\frac{\a^2}{N}} = O\brac{\frac{b \log^2 b}{n}}
\end{equation}
where we have used \eqref{eqn:xasymptotics} and \eqref{eqn:nmidef}.
It remains to show that \eqref{eqn:totalbigi} is negligible compared to \eqref{eqn:smalli2}. Suppose first that $1 \le b \le \sqrt{n}$. Then we have $ N =\Omega(\sqrt{n})$, $x =O(\log^3 n)$ and $\a \ge 10$. In this case \eqref{eqn:totalbigi} is
\[
     O\bigg( x (b+1) \brac{\frac{2x}{N}}^\a \bigg)= O\bigg( \log^3 n \cdot \sqrt{n} \brac{\frac{\log^3 n}{\sqrt{n}}}^{10} \bigg)
\]
which is negligible compared to \eqref{eqn:smalli2}. Now suppose $\sqrt{n} \le b \le n/\log^4 n$. Then $N = \Omega( \log^4 n)$, $x =O(\log^3 n)$ and $\a \ge 5 \log n$. In this case \eqref{eqn:totalbigi} is
\[
     O\bigg( x (b+1) \brac{\frac{2x}{N}}^\a \bigg)= O\bigg( \log^3 n \cdot \frac{n}{\log^4 n} \brac{\frac{\log^3 n}{\log^4 n}}^{5 \log n} \bigg).
\]
which is also negligible compared to \eqref{eqn:smalli2}. Thus the expected cost paid by Maker is $O\brac{\frac{b \log^2 b}{n}}$, completing the proof of Theorem \ref{thm:item} in the case where $2 \le b \le n / \log^4 n$. 

\subsection{Upper bound for $b=1$}
Finally we do the upper bound for $b=1$. Say that the cost $c_i$ of vertex $i$ is uniformly distributed on $[0,1]$. Breaker examines the vertices $1, \ldots, n$ in order and removes the first vertex $j$ such that $c_j \le b_j$ (for some deterministic numbers $b_1, \ldots b_n$). Let the random variable $B$ be this value of $j$. Maker then examines the vertices in order and takes the first vertex $i$ such that $i \neq B$ and $c_i \le m_i$ (for some deterministic $m_1, \ldots, m_n$). Let $M$ be this value of $i$.

The expected cost paid by Maker is precisely
\begin{align*}
    \E[c_{M}] & = \sum_{i=1}^{n} \Pr(M=i \; \wedge \; B>i) \E[c_i \mid  M=i \; \wedge \; B>i]\\
    & \qquad + \sum_{i=2}^n \sum_{j=1}^{i-1} \Pr(M=i \; \wedge \; B=j) \E[c_i \mid  M=i \; \wedge \; B=j].
\end{align*}
If Breaker plays optimally we will have $b_i\le m_i,i=1,2,\ldots,n$. Indeed, it would never be optimal for Breaker to take an item which Maker would not take if offered.  Now for all $i \ge 1$ we have
\[
\Pr(M=i \; \wedge \; B>i) =  \brac{\prod_{k=1}^{i-1}(1-m_k)} \cdot (m_i-b_i)
\]
(we regard an empty sum as 0 and an empty product as 1). We have
\[
\E[c_i \mid  M=i \; \wedge \; B>i] = \frac{m_i + b_i}{2}.
\]
Meanwhile we have for $1 \le j < i$ that 
\[
\Pr(M=i \; \wedge \; B=j) =  \brac{\prod_{k=1}^{j-1}(1-m_k)} \cdot b_j \cdot  \brac{\prod_{k=j+1}^{i-1}(1-m_k)} \cdot m_i,
\]
and
\[
\E[c_i \mid  M=i \; \wedge \; B=j] = \frac{m_i}{2}.
\]
Thus we let
\begin{align}
   c(\bf{b}, \bf{m}):= \E[c_{M}] & = \sum_{i=1}^n \brac{\prod_{k=1}^{i-1}(1-m_k)} \cdot \frac{m_i^2-b_i^2}{2} \nn \\
    & \qquad + \sum_{i=2}^n \sum_{j=1}^{i-1} \brac{\prod_{k=1}^{ j-1}(1-m_k)} \cdot b_j \cdot  \brac{\prod_{k=j+1}^{i-1}(1-m_k)} \cdot\frac{m_i^2}{2}\label{eqn:c}.
\end{align}

Let $\tm_i = \frac{2}{n-i+1}$ for $i \le n-1$, and $ \tm_n = 1$. We now try to maximize $c(\bf{b}, \widetilde{\bf{m}})$ over all choices for $\bf{b}$. 
\begin{claim}
    An optimal choice for $\bf{b}$ is given by
\begin{equation}\label{eqn:bistar}
  b_n=0, \quad b_{n-1}=\frac12, \quad \mbox{and for $i^*\le n-2$, } b_{i^*}  = \frac{2}{n-i^*-1} - \frac{2}{(n-i^*)(n-i^*-1)} \sum_{\ell=1}^{n-i^*} \frac{1}{\ell}.
\end{equation}
\end{claim}

\begin{proof}

We prove this claim by (backward) induction. We treat $b_n$ and $b_{n-1}$ as base cases. Note that $\tm_n = \tm_{n-1} = 1$ and $c(\bf{b}, \widetilde{\bf{m}})$ does not depend on $b_n$ at all and so we might as well fix $b_n=0$ which agrees with the claim. Now we turn to $b_{n-1}$. The only terms depending on $b_{n-1}$ in \eqref{eqn:c} are when $i=n-1$ for the first sum, and when $i=n, j=n-1$ for the second sum. Thus we obtain
\begin{align*}
    \frac{\partial c({\bf b}, {\bf \tm})}{\partial b_{n-1}} &=  -\brac{\prod_{k=1}^{n-2}(1-\tm_k)} \cdot b_{n-1}  + \brac{\prod_{k=1}^{n-2}(1-\tm_k)}     \cdot\frac{1}{2}\\ 
    &=  \left(\frac12 - b_{n-1}\right) \brac{\prod_{k=1}^{n-2}(1-\tm_k)} 
\end{align*}
and so $b_{n-1}= \frac12$ is an optimal choice, again agreeing with the claim. 

Suppose  for some $i^{**}$ with  $1 \le i^{**} \le n-2$ that  \eqref{eqn:bistar} gives the optimal value for $b_{i^*}$ for all $i^*>i^{**}$. We fix those values $b_{i^*}$ for all $i^*>i^{**}$ and determine the optimal value for $b_{i^{**}}$. 

Note that the product of $1-\tm_k$ is telescoping. For example when $i^* +1 \le i \le n$ we have 
\[
\prod_{k=i^{**}+1}^{ i-1}(1-\tm_k) = \prod_{k=i^{**}+1}^{ i-1} \frac{n-k-1}{n-k+1} = \frac{(n-i +1)(n-i)}{(n-i^{**})(n-i^{**}-1)}.
\]

Now
\begin{align*}
    \frac{\partial c({\bf b}, {\bf \tm})}{\partial b_{i^{**}}} &=  -\brac{\prod_{k=1}^{ i^{**}-1}(1-\tm_k)} \cdot b_{i^{**}} \\
    &  \qquad + \sum_{i=i^{**}+1}^{ n} \brac{\prod_{k=1 }^{ i^{**}-1}(1-\tm_k)} \cdot   \brac{\prod_{k=i^{**}+1}^{ i-1}(1-\tm_k)} \cdot\frac{\tm_i^2}{2}.
\end{align*}
So an optimal choice is 
\begin{align}
    b_{i^{**}} &= \sum_{i=i^{**}+1}^{ n}    \brac{\prod_{k=i^{**}+1}^{ i-1}(1-\tm_k)} \cdot\frac{\tm_i^2}{2} \nn\\
    &=\sum_{i=i^{**}+1}^{ n}   \frac{(n-i +1)(n-i)}{(n-i^{**})(n-i^{**}-1)} \cdot\frac{2}{(n-i+1)^2} \nn\\
    &= \frac{2}{(n-i^{**})(n-i^{**}-1)} \sum_{i=i^{**}+1}^{ n} \frac{n-i}{n-i+1}
\end{align}
which agrees with \eqref{eqn:bistar} proves the claim . 
\end{proof}

Supposing Maker plays by strategy ${\bf \tm}$ then the best Breaker can do is to play by the strategy given in \eqref{eqn:bistar}. We now estimate \eqref{eqn:c} for these strategies:
\begin{equation}\label{eqn:ceval1}
   \sum_{i=1}^n \brac{\prod_{k=1}^{i-1}(1-\tm_k)} \cdot \frac{ \tm_i^2-b_i^2}{2}  + \sum_{i=2}^n \sum_{j=1}^{i-1} \brac{\prod_{k=1}^{j-1}(1-\tm_k)} \cdot b_j \cdot  \brac{\prod_{k=j+1}^{i-1}(1-\tm_k)} \cdot\frac{\tm_i^2}{2} 
\end{equation}

When $i=n$, the term in the first sum is 0 since $\tm_{n-1}=1$. When $i=n$ for the second sum, the only nonzero term is when $j=n-1$ and so we get
\[
\brac{\prod_{k=1}^{n-2}(1-\tm_k)} \cdot \frac 12 \cdot \frac 12 =  \frac{3}{2n(n-1)} = O(n^{-2}).
\]
Thus \eqref{eqn:ceval1} becomes
\begin{align}
    &\sum_{i=1}^{n-1} \brac{\prod_{k=1}^{i-1}(1-\tm_k)} \cdot \frac{ \tm_i^2-b_i^2}{2}  + \sum_{i=2}^{n-1} \sum_{j=1}^{i-1} \brac{\prod_{k=1}^{j-1}(1-\tm_k)} \cdot b_j \cdot  \brac{\prod_{k=j+1}^{i-1}(1-\tm_k)} \cdot\frac{\tm_i^2}{2} +O(n^{-2})\nn\\
    &= \sum_{i=1}^{n-1} \frac{(n-i +1)(n-i)}{n(n-1)} \cdot \frac{ \tm_i^2-b_i^2}{2}  + \sum_{i=2}^{n-1} \sum_{j=1}^{i-1} \frac{(n-i +1)(n-i)(n-j+1)}{n(n-1)(n-j-1)} \cdot\frac{b_j \tm_i^2}{2} +O(n^{-2}).\label{eqn:ceval2}
\end{align}
We bound the first sum above. Note that 
\begin{equation}
  b_i, m_i = O\brac{(n-i)^{-1}}, \qquad   b_i = \tm_i + O\brac{(n-i)^{-2} \log(n-i)}
\end{equation}
and so
\[
\tm_i^2-b_i^2 = (\tm_i+b_i)(\tm_i-b_i) = O\brac{(n-i)^{-3} \log(n-i)}.
\]
Thus the first sum in \eqref{eqn:ceval2} is
\[
    \sum_{i=1}^{n-1} \frac{(n-i +1)(n-i)}{n(n-1)} \cdot \frac{ \tm_i^2-b_i^2}{2} =  O\brac{\sum_{i=1}^{n-1} \frac{\log(n-i)}{n^2(n-i)}} = O\brac{ \frac{\log^2n}{n^2}}.
\]
We turn to the second sum in \eqref{eqn:ceval2} which is
\begin{align}
&  \frac{1}{2n(n-1)} \sum_{i=2}^{n-1}  (n-i +1)(n-i)\tm_i^2\sum_{j=1}^{i-1}\brac{1+\frac{2}{n-j-1}} b_j    \nn\\
&= \frac{1}{2n(n-1)} \sum_{i=2}^{n-1} { \brac{4-\frac{1}{n-i+1}}}\ \sum_{j=1}^{i-1} \brac{\frac{2}{n-j} +O\bfrac{\log n}{(n-j)^{2}}}    \nn\\
& = \frac{1}{n(n-1)}\sum_{i=2}^{n-1} {\brac{4-\frac{1}{n-i+1}}}\brac{- \log\brac{\frac{n-i}{n}} + O\bfrac{\log n}{n-i}}\nn\\
& \sim \frac 4n
\end{align}
since $\sum_{i=2}^{n-1}  \log\brac{\frac{n-i}{n}} \sim n \int_0^1 \log(1-x) \; dx =-n $.

\section{Purchasing a $k$-clique}

In this section we prove Theorem \ref{thm:clique}. We start with the game \textsc{$3$-clique} without any phase restriction. Here and several times in the future we will use the Chernoff--Hoeffding  bound (see, for example, Theorem 23.6 in \cite{book})

\begin{theorem}[Chernoff--Hoeffding bound] \label{thm:chernoff}
    Let $X$ be distributed as $Bin(n, p)$ and $0 < \e < 1$. Then 
    \[
    \Pr(|X-np|>\e np) \le 2 \exp(-\e^2 np/3)
    \]
\end{theorem}

\subsection{Unrestricted \textsc{$3$-clique}}

Maker will first build a $n^{1/3}$-star by the time Maker has seen half of the edges. Maker can do this by arbitrarily choosing a root $v$ at the beginning, and then taking every offered edge incident with $v$ of cost at most $8(b+1)n^{-2/3}$.  Indeed, by Chernoff-Hoeffding \ref{thm:chernoff} we see that with failure probability at most $\exp(-\Omega(n^{1/3}))$ the number of such edges among the first $\binom n2 /2$ edges is at least half its expectation, i.e. $2(b+1)n^{1/3}$. Of course Breaker can take some, but Maker will get at least a $1/(b+1)$ fraction of them so Maker easily gets $n^{1/3}$ of them, completing the $n^{1/3}$-star rooted at $v$. The total cost paid by Maker for the star is at most $8(b+1)n^{-1/3}$. Maker still gets to look at the second half of the edges to complete the triangle. 

When Maker has finished building the $n^{1/3}$-star, Breaker has only taken $bn^{1/3}$ edges total. Unfortunately we have to deal with the possibility that Breaker has taken many edges in the neighborhood of $v$. Let $K$ be the set of leaves of our $n^{1/3}$-star. The expected number of edges  in $K$ among the second half of the edges in the random permutation of cost at most $40bn^{-1/3}\log n$ is $\frac12 \binom {n^{1/3}}2 \cdot \frac{40b\log n}{n^{1/3}} \sim 10bn^{1/3} \log n$. By Chernoff-Hoeffding \ref{thm:chernoff}, the probability that this number of edges is at most $5bn^{1/3} \log n$ (i.e. half its expectation) is at most $\exp\brac{-  10bn^{1/3} \log n / 12 }$ which is small enough to beat a union bound over choices of $K$, the number of choices here being $\binom n{n^{1/3}} \le \exp(n^{1/3} \log n)$. Thus, with high probability Maker is able to close the triangle using an edge of cost at most $20bn^{-1/3}\log n$.

Altogether, Maker pays a cost at most 
\[
8(b+1)n^{-1/3} + 20bn^{-1/3}\log n  \le 30 (b+1)n^{-1/3} \log n. 
\]

\subsection{$k$-phase restricted \textsc{$k$-Clique}}

In this version we consider the first $\frac1k \binom n2$ edges to be Phase 1, the next $\frac1k \binom n2$ edges to be Phase 2, and so on. If Breaker reaches the end of a phase before Maker, then Maker is allowed to take any edges Maker wants from the rest of that phase. Breaker is not allowed to begin the next phase until Maker has done this. 

In each of the first $k-3$ phases, Maker will build a $(k-3)$-clique with a large common neighborhood. Maker does this by building a star in each of these phases, each star contained inside the leaves of the last star. Finally in the last three phases Maker will build a triangle contained in the common neighborhood of the vertices in this $(k-3)$-clique, making a $k$-clique. 

\subsubsection{The first $k-3$ phases}

We describe how to build our stars. Set
\begin{equation}
  r:= n^{-\a_k }, \qquad  \ell_i := r \cdot \brac{\frac nr}^{2^{-i}}, 0 \le i \le k-3
\end{equation}
and note that for each $i$, $1 \le i \le k-3$ we have
\begin{equation}\label{eqn:identities}
    \frac{\ell_i^2}{\ell_{i-1}} = r= \ell_{k-3}^{-4/7}.
\end{equation}

 When we start Phase $i$ for $1 \le i \le k-3$, assume we have chosen a set $L_{i-1}$ of $\ell_{i-1}$ vertices (note that $\ell_0 = n$, so $L_0$ is all the vertices).  Maker chooses an arbitrary vertex $v_i \in L_{i-1}$, and during Phase $i$ Maker will take any edge of cost at most $10k(b+1)\ell_i/\ell_{i-1}$ with one end at $v_i$ and the other in $L_{i-1}$. Note that since $b = o\brac{ n^{11\a_k/4}}$ we have $10k(b+1)\ell_i/\ell_{i-1} = o(1)$. The probability that any particular edge appears in the first $i-1$ phases is $(i-1)/k \ge 1/k$. With high probability there are at least $(b+1) \ell_i$ such edges and so Maker gets $\ell_i$ of them at a total cost of at most $10k(b+1)\ell_i^2/\ell_{i-1} = 10k(b+1)r$ (using \eqref{eqn:identities}). For $1 \le i \le k-3$, $L_i$ will be the set of leaves of the star built at Phase $i$. 

 Altogether for Phases 1 through $k-3$, Maker pays at most $10(b+1)k^2 r$. If Maker builds a triangle in $L_{k-3}$ then Maker is done. 

\subsubsection{Phase $k-2$}

Set $\ell := \ell_{k-3}$, $L := L_{k-3}$. Maker's strategy now is to take any edge with both ends in $L$ whose cost is at most $10k(b+1)\ell^{-9/7}$ until Maker has $2 \ell^{5/7}$ such edges. The expected number of Phase $k-2$ edges that have cost at most $10k(b+1)\ell^{-9/7}$ is 
\[
\frac1k \binom {\ell}2 \cdot 10k(b+1)\ell^{-9/7} \sim 5(b+1)\ell^{5/7}
\]
and so an easy application of Chernoff-Hoeffding gives us that there are at least $2(b+1)\ell^{5/7}$ with exponentially small failure probability. Thus Maker will get at least $2\ell^{5/7}$ of them regardless of what Breaker does. Now we would like to claim that this set of $2\ell^{5/7}$ edges contains a matching on $\ell^{5/7}$ edges. Indeed, the number of vertices that are adjacent to two edges of cost at most $10(b+1)\ell^{-9/7}$ has expectation at most 
\[
\ell^3 \cdot \brac{10(b+1)\ell^{-9/7}}^2 = O\brac{\ell^{3/7}}
\]
and so we can easily remove $O\brac{\ell^{{ 4/7}}}$ of Maker's $2 \ell^{5/7}$ edges to get our matching. 

Maker ends Phase $k-2$ with a matching on $\ell^{5/7}$ edges. The total cost paid in Phase $k-2$ is at most 
\[
\sim 2 \ell^{5/7} \cdot 10k(b+1)\ell^{-9/7} = 20 (b+1)k\ell^{-4/7} = 20 (b+1)kr
\]
where we have used \eqref{eqn:identities}.

\subsubsection{Phase $k-1$}

In Phase $k-1$, Maker's strategy is as follows. Let $M$ be the matching on $m=\ell^{5/7}$ edges from Phase $k-2$. In Phase $k-1$ Maker will take any offered edge $e$ with the following properties, until Maker has $\ell^{4/7}$ of them:
\begin{enumerate}
    \item $e \subseteq L$ and is adjacent to some edge $e' \in M$,
    \item $e$ has cost at most $5(b+1)k^2\ell^{-8/7}$, and
    \item the edge $e''$ which would complete a triangle with $e, e'$ has not been offered to either player yet.
\end{enumerate}

By Chernoff it holds with exponentially small failure probability that for every edge $f \in M$ there are $\sim \frac {2}{k^2} \ell$ pairs of edges $f', f'' \subseteq L$ making a triangle with $f$ such that $f'$ is a Phase $k-1$ edge and $f''$ is a Phase $k$ edge. There are then 
\[
\sim \ell^{5/7} \cdot \frac {2}{k^2} \ell \cdot 5(b+1)k^2\ell^{-8/7} = 10(b+1)\ell^{4/7}
\]
such triples $f, f', f''$ where the cost of $f'$ is low enough that Maker would want it. In this scenario if  Maker is offered $f'$ during Phase $k-1$ then Maker would take it (since it is impossible that $f''$ has been offered to anyone yet). 

Technically, we will not condition on the event from the previous paragraph. We will condition on Phase $k-2$ completing as described, and then we will reveal the edges of Phase $k-1$ one by one as they are revealed by the players (i.e. whenever a player is offered an edge which has not been offered to either player yet). Of all the steps during which an edge is revealed which Maker would want, Maker gets at least a $(b+1)^{-1}$ fraction of them. For each edge $e'$ taken by Maker, conditional on having revealed all the edges up to $e'$ there is probability at least  $\frac 12$ that the edge $e''$ is a Phase $k$ edge. Indeed, the remaining unrevealed edges are a uniform random permutation of edges consisting of some of the phase $k-1$ edges and all of the phase $k$ edges. As discussed in the above paragraph, with exponentially small failure probability at least $\sim 10(b+1)\ell^{4/7}$ steps when such an edge $e'$ is revealed. By standard concentration arguments, w.h.p. we have that Maker gets at least 
\[
\sim 10(b+1)\ell^{4/7} \cdot (b+1)^{-1} \cdot \frac 12 > \ell^{4/7}
\]
edges $f'$ such that $f''$ is a Phase $k$ edge. The total cost paid in Phase $k-1$ is at most 
\[
\ell^{4/7} \cdot 5(b+1)k^2\ell^{-8/7} = 5(b+1)k^2r.
\]

\subsubsection{Phase $k$ and finishing the proof}

Now Maker just wants to take one of the $\ell^{4/7}$ possible edges which would close a triangle, and by Theorem \ref{thm:item} Maker can do that with an expected cost of $O( b \log ^2 b / \ell^{4/7} ) = O(b r \log ^2 b)$. The cost of all previous phases is absorbed into the big-O. This completes the proof of Theorem \ref{thm:clique}.

\section{Paths}
In this section we prove Theorem \ref{thm:paths}

We describe how to adapt the strategy of the second author and Pegden \cite{FP18}. They showed that for $b=0$ Maker can build a path connecting two given vertices with a total cost $n^{-2/3 + o(1)}$. Our strategy for Maker will essentially follow the strategy in \cite{FP18}, which still works even though Breaker will take some of the edges Maker wants. 

Let $k:= \log \log n$. Maker's strategy will be $3k$-phase restricted. In the first $k$ phases, Maker will build a sequence of trees $T_1\subseteq  \ldots \subseteq T_k$ containing $u$. Specifically, we let $T_0$ consist of just the vertex $u$, and for $1 \le i \le k$, at phase $i$ Maker will take any offered edge with one endpoint in $T_{i-1}$ whose cost is at most $(b+1)p$ where $p:=n^{-1+1/3k}$, until Maker has $((1-\e)np/3k)^i$ such edges, where $\e = n^{-1/9k}$. Note that since $b \le n^{1-1/k}$ we have that $(b+1)p\le 1$. Of course there is the possibility that Maker will not find this number of edges in which case phase $i$ fails. However, we will now see that is unlikely. Indeed, assuming phases $1, \ldots, i-1$ have succeeded, the number of vertices in $T_{i-1}$ satisfies
\[
((1-\e)np/3k)^{i-1} \le |V(T_{i-1})| = \sum_{j=0}^{i-1} ((1-\e)np/3k)^j \le (1+o(1)) ((1-\e)np/3k)^{i-1}.
\]
 The number of vertices $w$ outside of $T_{i-1}$ such that there is a phase $i$ edge from $w$ to $T_{i-1}$ of cost at most $p$ is distributed as $Bin(|V(T_{i-1})|(n- |V(T_{i-1})|), (b+1)p/3k)$. By Chernoff-Hoeffding, the probability that this random variable is less than $1-\e/2$ times its expectation $\m$ where
\[
\m = |V(T_{i-1})|(n- |V(T_{i-1})|) (b+1)p/3k \ge  ((1-\e)np/3k)^{i-1} (1-\e/3) (b+1) np/3k
\]
is at most $2\exp(-\e^2 \m / 12)=o(n^{-10})$. Thus with probability at least $1-o(n^{-10})$ there are at least $(1-\e/2)\m \ge (b+1)((1-\e)np/3k)^i$ edges that Maker would want in phase $i$. Since Maker will get at least a $1/(b+1)$ proportion of these edges, in this case phase $i$ is successful. Thus with probability at least $1-o(n^{-9})$ phases $1, \ldots, k$ are all successful. The amount paid by Maker during these phases is at most 
\[
|V(T_k)|(b+1)p = O\brac{((1-\e)np/3k)^{k} (b+1)p } = (b+1)n^{-2/3+o(1)}
\]

In phases $k+1, \ldots, 2k$, Maker builds a sequence of trees $T_1'\subseteq  \ldots \subseteq T_k'$ containing $v$. This works similarly the first $k$ phases with $v$ in place of $u$. Maker is able to build $T_k'$ containing  $((1-\e)np/3k)^{k}=n^{1/3+o(1)}$ vertices including $v$, and paying a total cost of at most $(b+1)n^{-2/3+o(1)}$. 

Now if $T_k, T_k'$ do not already intersect, Maker can very likely find a cheap edge connecting them in one of the phases $2k+1, \ldots 3k$ (we really do not even need to treat these as separate phases). Indeed, the expected number of edges from $T_k$ to $T_k'$ that cost at most $(b+1)\log^2 n / |V(T_k)| \cdot |V(T_k)'| = (b+1) n^{-2/3+o(1)}$ in the last $k$ phases is at least 
\[
|V(T_k)| \cdot |V(T_k)'| \cdot \frac 13 \cdot \frac{(b+1)\log^2 n}{|V(T_k)| \cdot |V(T_k)'|} = \frac 13 (b+1) \log^2 n.
\]
By Chernoff-Hoeffding there are easily $b+1$ such edges with high probability and so Maker gets one. Maker is done, and the total cost paid is $(b+1) n^{-2/3+o(1)}.$ This completes the proof of Theorem \ref{thm:paths}.

\section{Box game}
Here we prove Theorem \ref{thm:boxgame}. First assume the ordering is arbitrary, possibly adversarial and adaptive. Then Maker has a winning strategy if and only if $m \ge bn+1$. Indeed, if $m$ is that large then Maker can win by taking a ball if and only if it belongs to a box which has the fewest balls that do not yet belong to Breaker (here we say a ball belongs to Breaker if Breaker has taken it or if it has been rejected by Maker). After Breaker has taken $i$ turns there is never a box with more than $bi$ balls belonging to Breaker including those rejected by Maker by induction on $i$. Thus, Maker can always win if $m \ge bn+1$. Now suppose $m \le bn$. The ordering can be such that (at least towards the beginning) Breaker is always offered balls from the same box while Maker is only offered balls from the other $n-1$ boxes. Maker needs one ball from each of the other boxes so Maker is forced to use up at least $n-1$ turns. Meanwhile during these $n-1$ turns Breaker has accumulated $b(n-1)$ balls from the same box and wins the game on the next turn by taking $b$ more. 

Now assume the ordering of the balls is random, $b \ge b_0$, $\e >0$ is fixed and $m \le (1-\e)bn$. By Chernoff-Hoeffding, the probability that we do not see box 1 $b$ times over the course of the first $(1+\e/2)bn$ steps is at most $n^{-3}$. Likewise the probability that we see box 1 more than $b$ times or less than $(1-\e)b$ times in the first $(1-\e/2)bn$ steps has the same bound. This bound is small enough so that we can take the union bound over all boxes and all sequences of consecutive steps. Thus with high probability our sequence of balls is such that in any consecutive  $(1+\e/2)bn$ steps we see every ball $b$ times, and in any consecutive $(1-\e/2)bn$ steps we see every ball between $(1-\e)b$ and $b$ times.  

On Breaker's first turn, Breaker will always take $b$ balls from the first box. Breaker will only continue to take balls from the first box until Maker takes a ball from that box (which requires Maker to look ahead of Breaker in the ordering, since Breaker has taken all balls in box 1 before Breaker's pointer). If Maker does take a ball from the first box, then Breaker chooses any other box which Maker has no balls in yet and starts taking balls from that box until Maker gets one, and then Breaker chooses a new box to focus on. Every time Maker takes Breaker's box, this requires Maker's pointer to go ahead of Breaker's. Breaker only takes balls that are ahead of Maker's pointer so whenever Maker gets ahead, Breaker will skip over some balls to get ahead.  Breaker needs to look at at least $(1-\e/2)bn$ balls each turn in order to get $b$ from the desired box. In the process of looking at those $(1-\e/2)bn$ balls, Breaker must pass over at least $(1-\e)b$ balls in every other box. If Maker allows Breaker to take balls from the same box for $k \ge 1$ turns, then Breaker has looked at least $k(1-\e/2)bn$ balls during these $k$ turns. Maker cannot get ahead of Breaker to take that box without passing over $k(1-\e)b$ balls in every other box. Therefore when Breaker picks a new box to focus on, over these $k$ turns Maker will have rejected $k(1-\e)b$ balls from that box already. By induction, after a total of $j$ turns there is always a box which has at least $j(1-\e)b$ balls which are unavailable to Maker. This completes the proof of Theorem \ref{thm:boxgame}.

\section{Closing remarks}

Our results are probably not optimal. In Theorem \ref{thm:item} the $\log^2 b$ factor is likely not necessary. Indeed, it is natural to expect that the optimal strategy for Maker gives an expected cost of $\sim 2(b+1)/n$. This is because, as Frieze and Pegden \cite{FP18} showed, when $b=0$ the optimal expected cost is $\sim 2/n$. Heuristically, Breaker should be able to inflate that by a factor $b+1$ since Maker can only get one out of every $b+1$ cheap items. However, in Theorem \ref{thm:item} we were only able to prove it for $b=1$. 

Regarding Theorem \ref{thm:clique}, Frieze and Pegden \cite{FP18} showed that for $b=0$ the expected cost of a triangle is at least $n^{-4/7}$ and that the $-4/7$ power is optimal. Thus it is somewhat satisfying to know that our result gets the same expected cost inflated by the factor $b$. However while it is natural to assume that this $b$ factor is necessary, it would be nice to see a strategy for Breaker which matches our upper bound. There is also the analogous question about the $b$ factor in Theorem \ref{thm:paths}. Also, for Theorems \ref{thm:clique} and \ref{thm:paths} it is unclear whether the phase restrictions are strictly necessary. 

Finally, in Theorem \ref{thm:boxgame} for the randomly ordered \textsc{BoxGame}, we believe the condition $b \ge b_0$ is unnecessary and the result should hold for all $b \ge 1$. We use $b \ge b_0$ to establish that long-enough intervals in the ordering of balls are all nicely concentrated in terms of the number of times each ball appears. For smaller $b$ we would not have the necessary concentration for our argument to go through. It seems that Breaker should actually be able to take advantage of this lack of concentration, but we are unable to prove it at the moment. 

\bibliographystyle{abbrv}

\end{document}